%% file: fibered_ribbon_knot_with_right_veering_pseudo_Anosov_monodromy.tex
\documentclass[12pt]{amsart}
\usepackage[utf8]{inputenc}

\usepackage{amssymb, amsmath, amsthm}
\usepackage[margin=1in]{geometry}
\usepackage{tikz, tikz-cd}
\usepackage{graphicx}
\usepackage{enumerate}
\usepackage{comment}
\usepackage{bbold}
\usepackage{subcaption}
\usepackage{pinlabel}
\usepackage{hyperref} 
\usepackage{thmtools}
\usepackage{thm-restate}
\usepackage{bm}
\usepackage{float}
\hypersetup{
    colorlinks = false, 
    linktoc = all,     
    linkcolor = darkgray,  
}
\usepackage[all,cmtip]{xy}
\usepackage{cleveref}

\usepackage[bordercolor=gray!20,backgroundcolor=blue!10,linecolor=blue!50,textsize=footnotesize,textwidth=0.8in]{todonotes}
\usepackage{pinlabel}

\input{commands.tex}

\title{hyperbolic fibered slice knots with right-veering monodromy}

\author[D. He]{Dongtai He}
\address {Shanghai Center for Mathematical Science, Fudan University, Shanghai, China 200438}
\email{dongtaihe@fudan.edu.cn}

\begin{document}

\maketitle
\begin{abstract}
We construct a hyperbolic fibered slice knot with right-veering monodromy, which disproves a conjecture posed in \cite{hubbard2021braids}.
\end{abstract}
\section{introduction}

We give a negative answer to Question 8.2 posed in \cite{hubbard2021braids}.

\begin{question}
If K is a hyperbolic fibered slice knot, does the fractional Dehn twist coefficient (FDTC) of the monodromy vanish?
\end{question}

We construct a hyperbolic slice fibered knot $K'$ with positive FDTC. 

\subsection{Motivation.} The authors \cite{hubbard2021braids} observes that many low-crossing slice fibered knots have zero FDTC. Any $(p,1)-$cable of slice fibered knot is still slice fibered, whereas FDTC of $(p,1)-$cable equals to $1/p$. The authors therefore ask the above question about hyperbolic fibered slice knots.  \\

Baldwin, Ni and Sivek \cite[Corollary 1.7]{baldwin2022floer} prove the following related proposition in terms of the $\tau-$invariant in Heegaard Floer homology:

\begin{proposition}\label{thin}
If $K\subset S^{3}$ is a fibered knot with thin knot Floer homology such that $\tau (K)<g(K)$, then FDTC vanishes.
\end{proposition}

The $\tau-$invariant vanishes for slice knots. Proposition \ref{thin} explains the case for low-crossing fibered slice knots because many of those are thin.\\

We have an immediate corollary:

\begin{corollary}
The knot Floer homology of $K'$ is not thin.
\end{corollary}

\subsection{Organization.} We follow the recipe of Kazez and Roberts \cite{Kazez-Roberts} to construct hyperbolic fibered knots with positive FDTC. The search for ribbon knot is inspired by the work of Hitt and Silver \cite{hitt1991ribbon}. In section \ref{2} we review Nielson-Thurston classification of surface automorphism and examples from Kazez and Roberts. We construct our example $K'$ in section \ref{ribbon}.

\section{monodromies of fibered knots in $S^{3}$}\label{2}

\subsection{Surface automorphism.} We first recall Nielson-Thurston classfication of surface automorphisms:

\begin{theorem}
\cite{casson1988automorphisms, thurston1988geometry} Let $S$ be an oriented hyperbolic surface with geodesic boundary, and let $h\in Aut(S, \partial S)$. Then $h$ is freely isotopic to either
\begin{enumerate}
\item a pseudo-Anosov homeomorphism $\phi$ that preserves a pair of geodesic laminations $\lambda^{s}$ and $\lambda^{u}$,
\item a periodic homeomorphism $\phi$ such that $\phi^{n}=id$ for some $n$,
\item a reducible homeomorphism $h'$ that preserves a maximal collection of simple closed geodesic curves in $S$. To avoid overlap, we consider $h$ reducible only when it is not periodic.  
\end{enumerate}
\end{theorem}

In particular, we only regard $h$ as reducible only if it is not periodic to avoid overlap. Let $\Phi:S\times [0,1]\rightarrow S$ be an isotopy from $h$ to its Thurston representative $\phi$.  Considering the restriction of $\Phi$ to the boundary $\partial S$, we have a homeomorphism:
$$\Phi_{\partial}:\partial S\times [0,1]\rightarrow\partial S\times [0,1]$$ defined by $\Phi_{\partial}(x,t)=(\Phi_{t}(x),t)$. The fractional Dehn twist coefficient $c(h)$ can be defined as the winding number of the arc $\Phi_{\partial}(\theta\times [0,1])$. Nielson-Thurston classification guarantees that $c(h)\in\mathbb{Q}$. \\

Thurston proved that a fibered knot is hyperbolic if and only if its monodromy is pesudo-Anosov. Fractional Dehn twist coefficient is closely related to the following notion of right-veeringness.

\begin{definition}
\cite{HKM07} A homeomorphism $h\in Aut(S,\partial S)$ is called right-veering if for every based point $x\in\partial S$ and every properly embedded arc $\alpha$ starting at $x$, $h(\alpha)$ is to the right of $\alpha$, after isotoping $h(\alpha)$ so that it intersects $\alpha$ minimally. Similarly, $h$ is called left-veering if $h(\alpha)$ is to the left of $\alpha$.
\end{definition}

\begin{proposition}
\cite{HKM07} $h$ is right-veering if and only if $c(h)>0$ for every component of $\partial S$, and $h$ is left-veering if and only if $c(h)<0$ for every component of $\partial S$.
\end{proposition}

If $c(h) = 0$, one can find two arcs such that one is moved by $h$ to the right and the other to the left. The significance of right-veeringness is highlighted by the following theorem of Honda, Kazez and Matić:

\begin{theorem}
\cite{HKM09} Every open book that is compatible with a tight contact structure is right-veering.
\end{theorem}

A large source of examples of reducible right-veering homeomorphism comes from the class of fibered cable knots. Indeed, if $h$ is the monodromy of a fibered $(p,q)-$cable knot $K_{p,q}$ with Seifert surface $S$, then $c(h)=1/pq$ and $h$ is reducible. Let $\{C_{i}\}$ be the collection of curves preserved by $h'$. $\{C_{i}\}$ partitions $S$ into subsurfaces $\{S_{j}\}$ permuted by $h'$. Let $S_{0}$ be the subsurface containing $\partial S=K_{p,q}$, then $h'|_{S_{0}}$ is periodic. Kazez and Roberts characterize the monodromy $h$ of a fibered knot $K$ in $S^{3}$ in the following theorem:

\begin{theorem}
\cite{Kazez-Roberts}
\begin{enumerate}
\item If $h$ is periodic, then $K$ is the unkont or a $(p,q)-$torus knot.
\item If $h$ has a reducible Thurstion representative $h'$ with periodic $h'|_{S_{0}}$, then $K$ is a $(p,q)-$cable knot, and $c(h)=1/pq$.
\item \cite{gabai1997problems} If $h$ is either pseudo-Anosov or reducible with $h'|_{S_{0}}$ pseudo-Anosov. Then either $c(h)=0$ or $c(h)=1/r$, where $2\leq |r|\leq 4g(K)-2$.
\end{enumerate}
\end{theorem}

\begin{corollary}
$c(h)=0$ or $1/r$ for some integer $r$, $|r|\geq 2$. In particular, $|c(h)|\leq 1/2$.
\end{corollary}

In particular, the $(2,1)-$cable of a fibered knot in $S^{3}$ has its monodromy attaining maximum FDTC. We review hyperbolic case in the next section.

\subsection{Stallings' twist and (2,1)-cable.} 
Let $U$ be an unknot properly embedded in a surface $F$. We say $U$ is untwisted relative to $F$ if $U$ bounds a disk transverse to $F$ along $U$. A Stallings' twist \cite{Stallings} is a surgery along such an untwisted $U$. Kazez and Roberts apply Stallings' twist on $(2,1)-$cables to produce hyperbolic fibered knots with maximum FDTC $= 1/2$.\\

Let $(S,h)$ be an open book decomposition of $S^{3}$ with connected binding $K$, where $h$ is pseudo-Anosov and $c(h)=0$. Let $K_{2,1}$ be the $(2,1)-$cable of $K$. The fibered surface $\Sigma$ of $K_{2,1}$ can be viewed as the union of two copies $S_{0}$, $S_{1}$ of $S$ connected by a 1-handle. Let $H$ be the monodromy of this new open book. \\

We choose a simple closed curve $C$ in $\Sigma$ such that $C_{0}=C\cap S_{0}$ and $C_{1}=C\cap S_{1}$ are two essential arcs. Moreover, we require $C_{i}$ to be nonseparating in $S_{i}$. Let $T_{C}$ be the right-handed Dehn twist along $C$ and $H' = T_{c}\circ H$.

\begin{theorem}
\cite{Kazez-Roberts} $H'$ is pseudo-Anosov and $c(H')=1/2$.
\end{theorem}

\section{ribbon fibered knot} \label{ribbon}
\begin{figure}[H]
\centering
\includegraphics[scale=0.5]{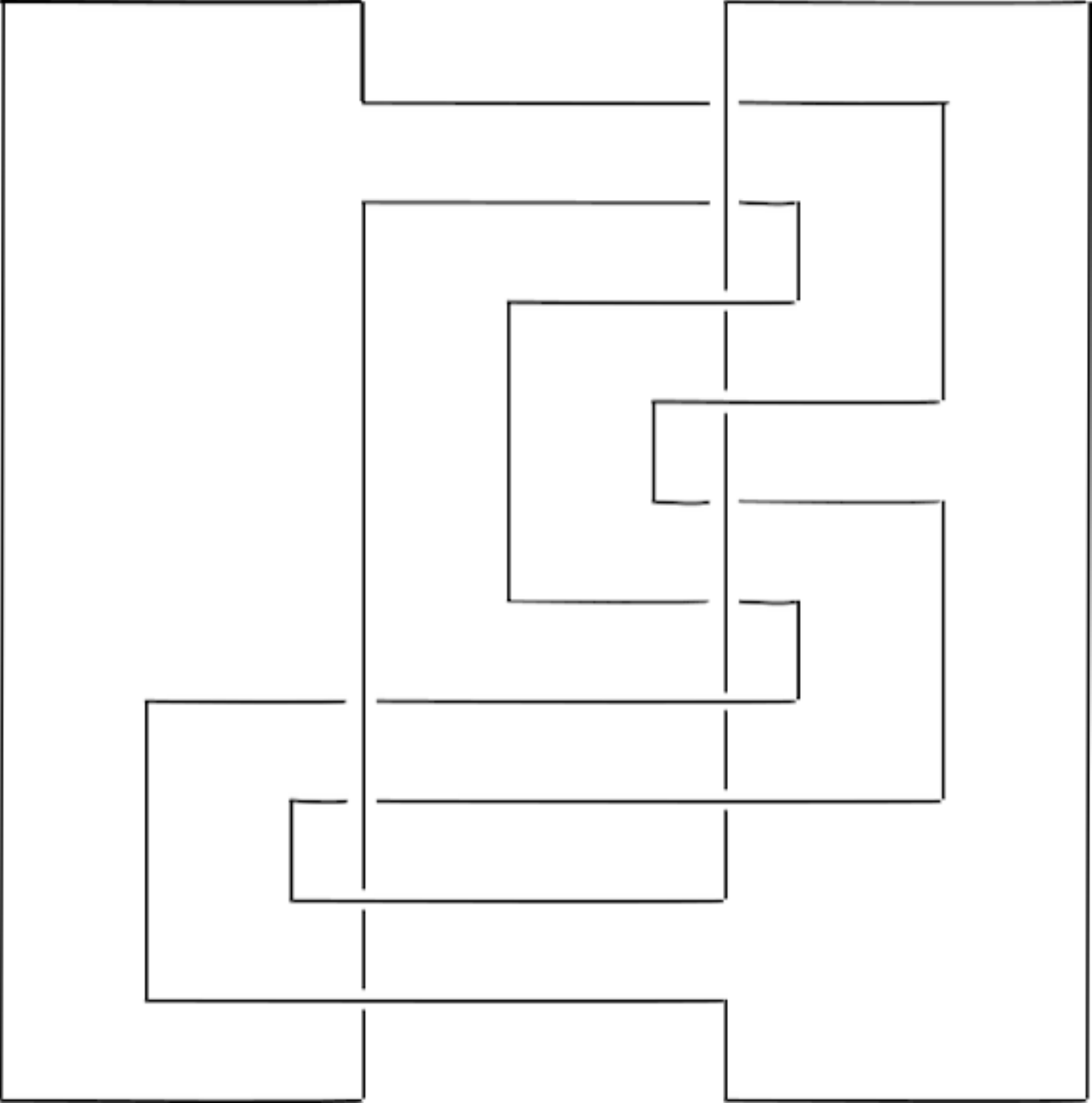}
\caption{A ribbon diagram for the knot $K=10_{153}$}\label{10_153}
\end{figure}

We are ready to construct a hyperbolic ribbon fibered knot with positive FDTC. Let $K$ be the knot $10_{153}$ from Rolfsen's knot table. $K$ is a hyperbolic ribbon fibered knot with 3-genus 3. Figure \ref{10_153} is a ribbon diagram for $10_{153}$. 

\begin{figure}[H]
\centering
\includegraphics[scale=0.5]{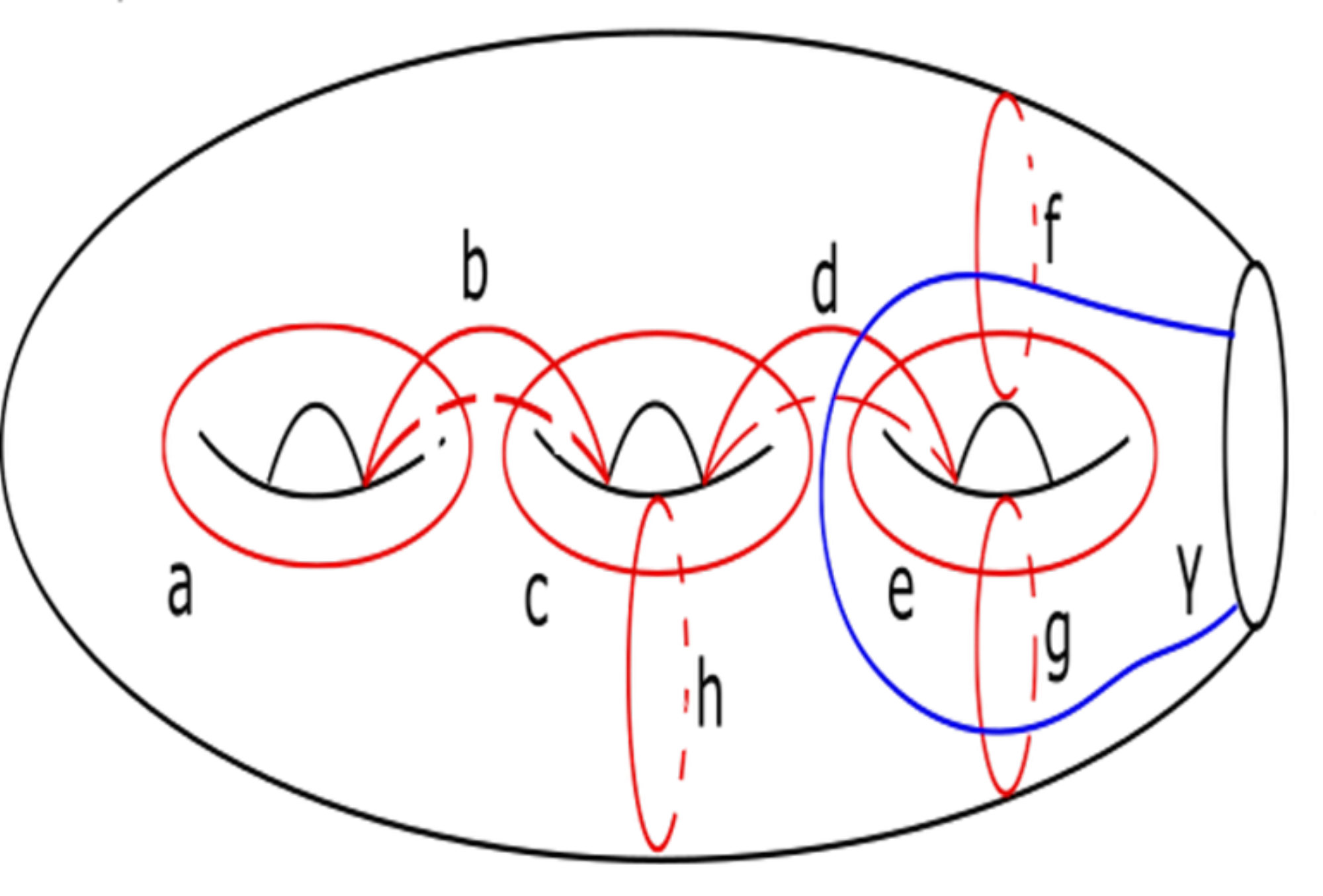}
\caption{Monodromy of the fibered knot $10_{153}$. $h$ can be presented as the word $\mathbf{abcBEGhcd}$, where $x$ denotes a right-handed Dehn twist about $x$ and $X$ denotes a left-handed Dehn twist about $x$. A word is read from right to left so that $aB$ means perform a left-handed Dehn twist about $b$ then perform a right-handed Dehn twist about $a$.}\label{monodromy}
\end{figure}

Let $h$ denote the monodromy. According to \cite{KnotInfo}, $h$ can be presented as decribed in Figure \ref{monodromy}. One can see that $h$ is neither right-veering nor left-veering by choosing different endpoints of $\gamma$. Therefore, $c(h)=0$.

\begin{figure}[H]
\centering
\includegraphics[scale=0.5]{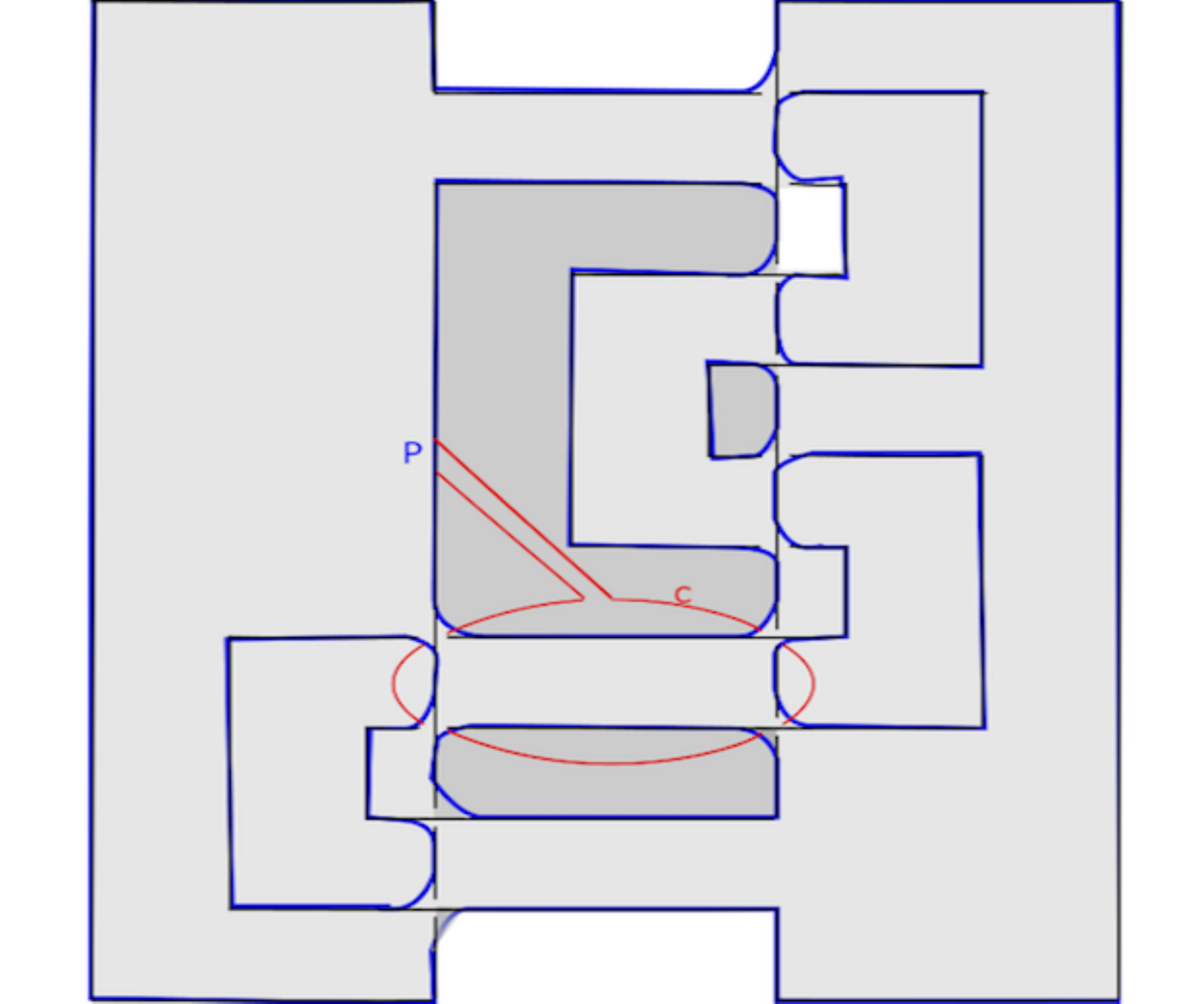}
\caption{The surface obtained by Seifert's algorithm has genus 3. $c$ is a non-seperating properly embedded arc on the surface.}\label{stallings}
\end{figure}

A Seifert surface $S$ of $K$ can be obtained by Seifert's algorithm as explained in Figure \ref{stallings}. The genus of $F$ is 3 so that $F$ is the fibered surface.\\ 

\begin{figure}[H]
\centering
\includegraphics[scale=0.5]{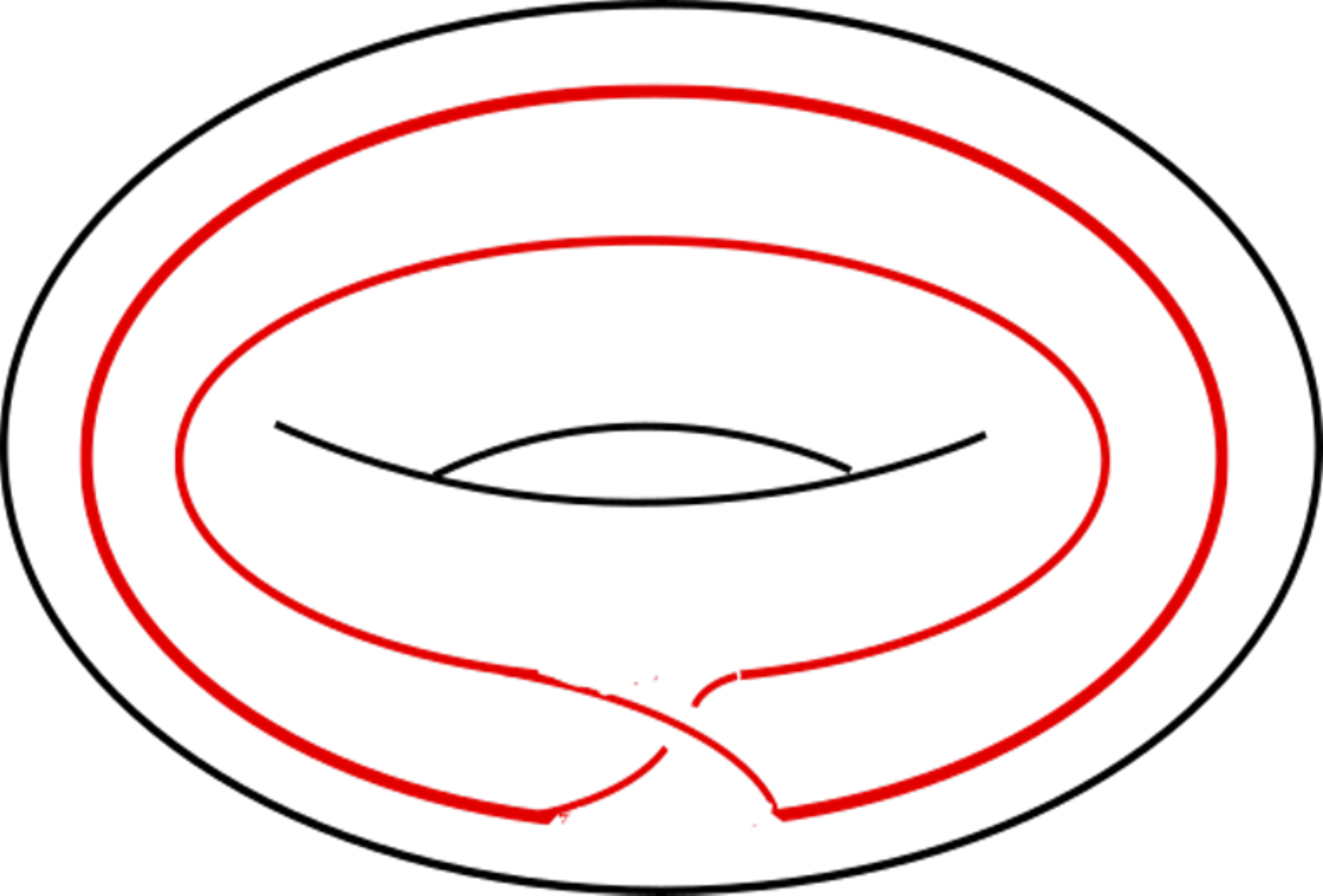}
\caption{$(2,1)-$cable of $K$.}\label{cable}
\end{figure}

Let $K_{2,1}$ be the $(2,1)-$cable of $K$ (Figure \ref{cable}). The twisted band connecting the two copies of $K$ is added at $p$. $c$ is a nonseperating properly embedded arc on the fibered surface of $K$. $K_{2,1}$ is also fibered whose fibered surface $\Sigma$ can be obtained by connecting two copies of $S$ with the same twisted band at $p$. Then define a simple closed curve $C$ to be a band sum of the two copies of $c$ along an arc running across the twisted band. 

Let $T_{C}$ denote the right-handed Dehn twist along $C$, and denote the resulting fibered knot $K'$. By \cite[Corollary 4.6]{Kazez-Roberts}, the monodromy $T_{C}\circ H$ is pseudo-Anosov and right-veering with $c(T_{C}\circ H)=\frac{1}{2}$. \\

\begin{figure}[H]
\centering
\includegraphics[scale=0.5]{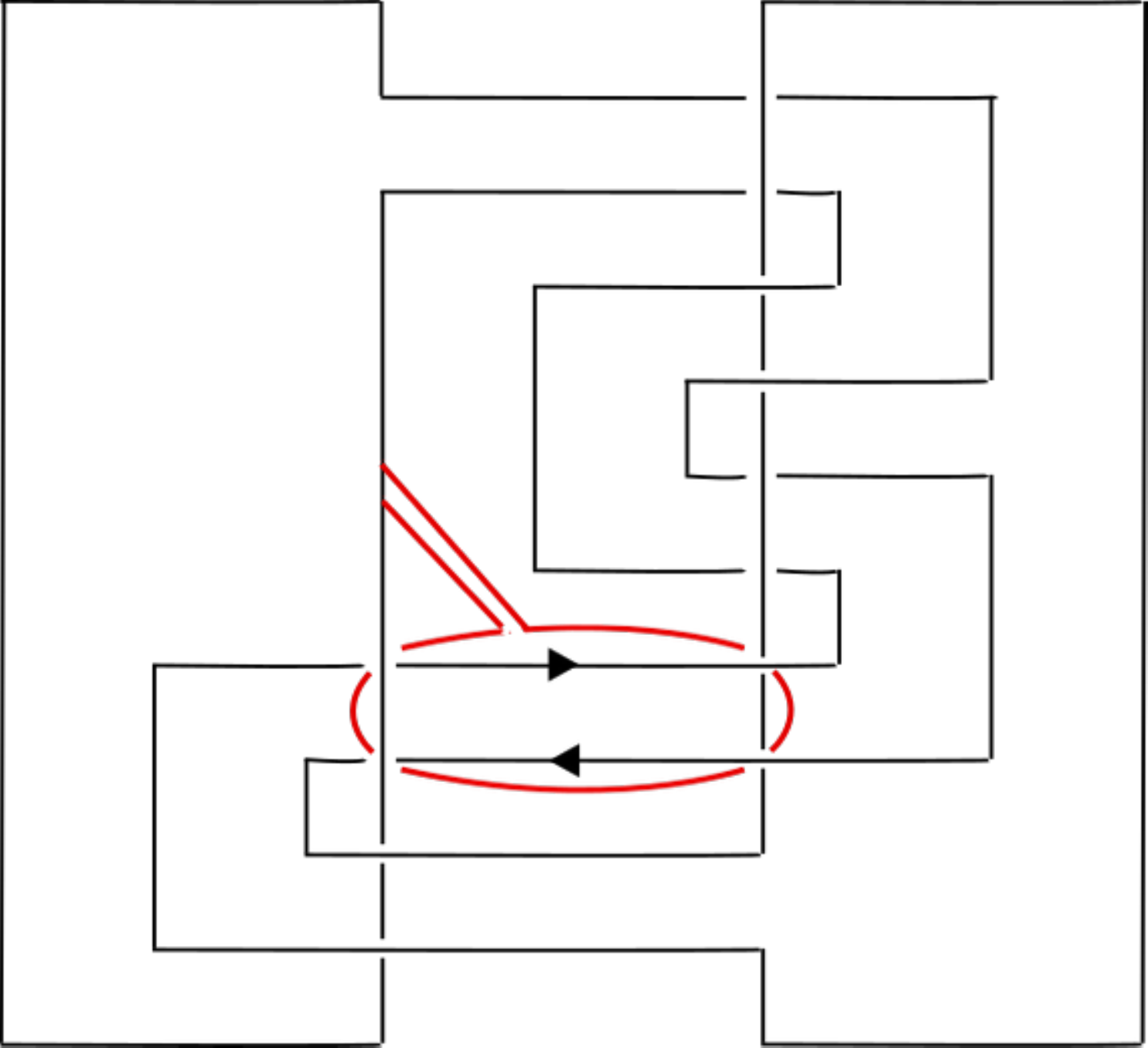}
\caption{The curve $C$ winds around two copies of a ribbon band. This figure shows one of the ribbon. The other ribbon is on the other copy from the $(2,1)-$cable}.\label{band}
\end{figure}

Recall that $K=10_{153}$ is a ribbon knot, so is $K_{2,1}$. $C$ is an unknotted untwisted curve. Performing a right-handed Dehn twist along $C$ has the same effect on $(S^{3}, K_{2,1})$ applying a $(-1)-$surgery along $C$. The resulting manifold is still $S^{3}$ and we have a new knot $K'$. $C$ winds around two copies of a ribbon band (Figure \ref{band}).

\begin{figure}[H]
\centering
\includegraphics[scale=0.5]{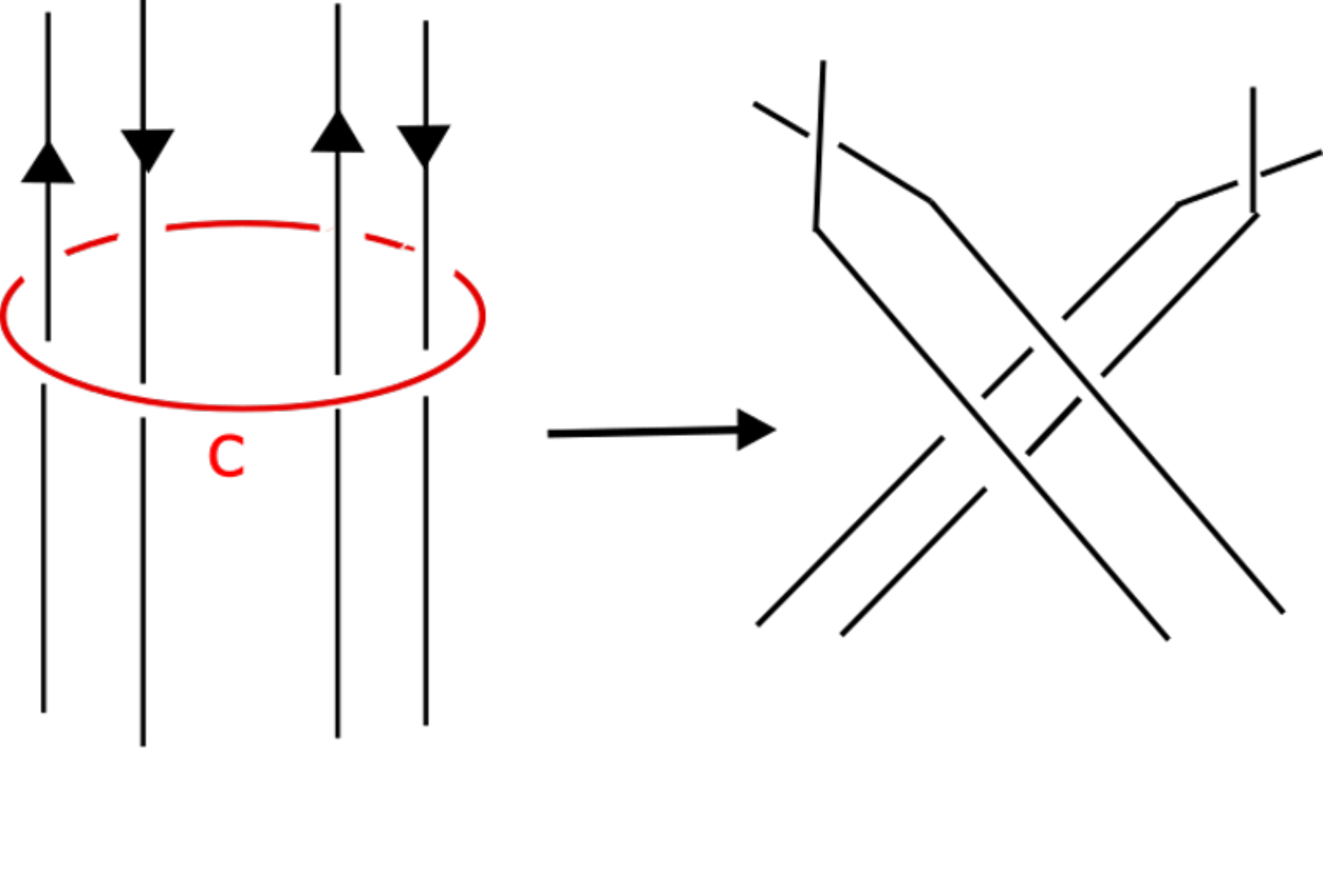}
\caption{The effect of $(-1)-$ surgery along $C$ after isotopy.}\label{surgery}
\end{figure}

Figure \ref{surgery} illustrates the effect of $(-1)-$surgery along $C$ to the ribbon bands. The resulting knot $K'$ is still a ribbon knot.
 
\begin{theorem}
$K'$ is a hyperbolic ribbon fibered knot with FDTC$=1/2$. Hence, the monodromy is right-veering.
\end{theorem} 

\bibliographystyle{alpha}
\bibliography{bib}
\end{document}

%% file: commands.tex
\theoremstyle{definition}
\newtheorem{theorem}{{Theorem}}[section]

\newtheorem{proposition}[theorem]{{Proposition}}
\newtheorem{definition}[theorem]{{Definition}}
\newtheorem{corollary}[theorem]{{Corollary}}

\newtheorem{question}[theorem]{Question}
